\documentclass[12pt,reqno,a4wide]{amsart}

\usepackage{tikz} 
\usepackage{calrsfs}
\usepackage[charter]{mathdesign}
\usepackage{hyperref}

\usetikzlibrary{decorations.pathreplacing}
\usetikzlibrary{fit}

\usepackage{pgfplots}

\allowdisplaybreaks

\DeclareMathOperator{\Res}{Res}

\begin{document}
\bibliographystyle{plain}

%
%

	\title[Convolution identities for {T}ribonacci numbers]
	{Convolution identities for {T}ribonacci numbers via the diagonal of a bivariate generating function}

	\author[H. Prodinger ]{Helmut Prodinger }
	\address{Department of Mathematics, University of Stellenbosch 7602, Stellenbosch, South Africa}
	\email{hproding@sun.ac.za}

	
	\begin{abstract}
Convolutions  for {T}ribonacci numbers involving binomial coefficients are treated with ordinary generating functions and the diagonalization method of Hautus and Klarner. In this way, the relevant generating function can be established, which is rational. The coefficients can also be expressed. It is sketched how to extend this to Tetranacci numbers and similar quantities.
	\end{abstract}
	

\maketitle


%
%

\section{Introduction}

Komatsu~\cite{Komatsu} treats
\begin{equation*}
\sum_{k=0}^n\binom nk T_kT_{n-k},
\end{equation*}
with Tribonacci numbers $T_m$, using exponential generating functions. This sequence, however has a rational generating function, which one does not see from this treatment. 

We present a method, based on ordinary generating functions, and the diagonal method of Hautus and Klarner~\cite{HK}, that provides this rational generating function.

As a warm-up, we will discuss the analogous question with Fibonacci numbers first.

\section{A warm-up}
We consider the double sequence
\begin{equation*}
	\sum_{k=0}^n\binom nk F_kF_{m-k},
\end{equation*}
with Fibonacci numbers $F_i$, and eventually specialize to $m=n$ (hence the name diagonal method).
We set up and compute a double generating function:
\begin{align*}
H(x,y)&=\sum_{n\ge0}\sum_{k=0}^n\binom nk F_k\sum_{m\ge k}y^mF_{m-k}\\
&=\frac{y}{1-y-y^2}\sum_{n\ge0}\sum_{k=0}^n\binom nk F_ky^k\\
&=\frac{y}{1-y-y^2}\sum_{k\ge0}\frac{x^k}{(1-y)^{k+1}} F_ky^k\\
&=\frac{xy^2}{(1-2x+x^2-xy-x^2y+x^2y^2)(1-y-y^2)}.
\end{align*}
Now, in order to pull out the diagonal, following Hautus and Klarner, we consider
\begin{equation*}
H(xt,\tfrac1t)\tfrac1t=\frac{z}{(1-2zt+z^2t^2-z+tz^2-z^2)(t^2-t-1)}
\end{equation*}
and look at the poles
\begin{equation*}
\frac{1\pm\sqrt5}{2}; \qquad \frac1z-\frac{1\pm\sqrt5}{2}.
\end{equation*}
One needs to consider the residues, but  the third and fourth poles go to infinity for $z\to\infty$.
So using Cauchy's integral theorem, these poles lie outside and do not need to be considered.

We find
\begin{equation*}
\Res(H(xt,\tfrac1t)\tfrac1t ;t=\tfrac{1+\sqrt5}{2})+\Res(H(xt,\tfrac1t)\tfrac1t ;t=\tfrac{1-\sqrt5}{2})=
\frac{z^2}{(1-z)(1-2z-4z^2)},
\end{equation*}
which is the rational function of interest. We can do a little bit more. Writing
\begin{equation*}
\alpha;\beta=\frac{1\pm\sqrt5}{2}
\end{equation*}
as usual, we can decompose
\begin{equation*}
\frac{z^2}{(1-z)(1-2z-4z^2)}=-\frac25\frac1{1-z}+\frac15\frac{1}{1-2\alpha z}+\frac15\frac{1}{1-2\beta z}.
\end{equation*}
From this we conclude that
\begin{equation*}
	\sum_{k=0}^n\binom nk F_kF_{n-k}=-\frac25+\frac152^nL_n,
\end{equation*}
with Lucas numbers $L_m=\alpha^m+\beta^m$.

\section{Tribonacci numbers}

Now consider
\begin{equation*}
	\sum_{k=0}^n\binom nk T_kT_{n-k},
\end{equation*}
where the Tribonacci numbers are given via the generating function
\begin{equation*}
\sum_{n\ge0}T_nz^n=\frac{1}{1-z-z^2-z^3}.
\end{equation*}
Exactly as in the warm-up section, we find the bivariate generating function
\begin{equation*}
G(x,y)=\frac{y}{1-y-y^2-y^3}\frac1{1-x}\frac{w}{1-w-w^2-w^3}\Big|_{w=\frac{xy}{1-x}}.
\end{equation*}
We look at $G(zt,\tfrac 1t)\tfrac 1t$ and its poles. There are 6 poles, but only 3 of them need to be considered, namely the roots
of $1-t-t^2-t^3$, call them $s_1, s_2, s_3$. While they are not particulary appealing, the computer can handle them. After heavy simplifications with a computer, we get
\begin{align*}
\Res(G(zt,\tfrac1t)\tfrac1t;t=s_1)&+\Res(G(zt,\tfrac1t)\tfrac1t;t=s_2)+\Res(G(zt,\tfrac1t)\tfrac1t;t=s_3)\\
&=\frac{1}{11}\frac{1+z+10z^2}{1-2z-4z^2-8z^3}-\frac{1}{11}\frac{1+z-8z^2}{1-2z+2z^3}.
\end{align*}
The coefficients of the first term are easy:
\begin{equation*}
\frac1{11}\Bigl(2^{n+1}T_{n+1}+\frac122^nT_n+\frac522^{n-1}T_{n-1}\Bigr).
\end{equation*}
For the second term we note that
\begin{equation*}
\frac{z^3}{1-z-z^2-z^3}\Big|_{z=\frac{-x}{1-x}}=\sum_{n\ge2}T_{n-2}\Big(\frac{-x}{1-x}\Big)^n=\frac{-x^3}{1-2x+2x^3}.
\end{equation*}
Setting
\begin{equation*}
\frac{1}{1-2z+2z^3}=\sum_{n\ge0}U_nz^n,
\end{equation*}
we find that
\begin{align*}
U_{m-3}&=-[x^m]\sum_{n\ge2}T_{n-2}\Big(\frac{-x}{1-x}\Big)^n=-[x^{m-n}]\sum_{n\ge2}T_{n-2}(-1)^n\Big(\frac{1}{1-x}\Big)^n\\
&=-[x^{m-n}]\sum_{n\ge2}T_{n-2}(-1)^n\sum_{k\ge0}\binom{n+k-1}{k}x^k\\
&=\sum_{n\ge2}T_{n-2}(-1)^{n-1}\binom{m-1}{n-1},
\end{align*}
or, nicer
\begin{align*}
	U_{m}=\sum_{k\ge1}T_{k-1}(-1)^{k}\binom{m+2}{k}.
\end{align*}
Then we can write the second contribution as
\begin{equation*}
-\frac{1}{11}[z^n]\frac{1+z-8z^2}{1-2z+2z^3}=-\frac{1}{11}(U_n+U_{n-1}-8U_{n-2}).
\end{equation*}

Komatsu~\cite{Komatsu} also considered arbitrary initial conditions for the Tribonacci numbers. That 
can be done with generating functions as well; only the numerator changes, and the 3 relevant poles stay the
same. The reader is invited to do a few experiments herself.

\section{Tetranacci and more}

Without giving details of the computations, the rational generating functions for Tetranacci number
(generating function $z/(1-z-z^2-z^3-z^4)$)is
\begin{equation*}
 {\frac {2{z}^{2} \left( -{z}^{3}-2 {z}^{4}+8 {z}^{5}+6 {z}^{6}+4
		 {z}^{7}+1-2 z-{z}^{2} \right) }{ \left( 16 {z}^{4}+8 {z}^{3}+4 {
			z}^{2}+2 z-1 \right)  \left( {z}^{6}+6 {z}^{5}-4 {z}^{4}-3 {z}^{3}
		-{z}^{2}+3 z-1 \right) }}.
\end{equation*}

The next instance is
\begin{equation*}
 {\frac {-2{z}^{2} \left( -{z}^{3}-{z}^{4}-25 {z}^{6}+19 {z}^{8}+52
		 {z}^{10}+40 {z}^{9}-1+3 z \right) }{ \left( 32 {z}^{5}+16 {z}^{4
		}+8 {z}^{3}+4 {z}^{2}+2 z-1 \right)  \left( 4 {z}^{10}-4 {z}^{9}-
		15 {z}^{8}-12 {z}^{7}+25 {z}^{6}-2 {z}^{4}-4 {z}^{3}-3 {z}^{2}+4
		 z-1 \right) }},
\end{equation*}
but after that the computations become too heavy to be reported here.

\bibliographystyle{plain}

\end{document}